\documentclass{article}

\usepackage{amsmath,amssymb,amsfonts}
\usepackage{url}
\newcommand{\qed}{\hfill \ensuremath{\Box}}
\DeclareMathOperator*{\ssum}{\sum\cdots\sum}

\begin{document}

\title{    More on the identity of Chaundy and Bullard	}
\author{		D. Aharonov and U. Elias			}
\date{}  
\maketitle


In  \cite{CB}  Chaundy and Bullard proved the identity
\begin{equation}							\label{eq:CB}							
	1 = x^{k+1} \sum_{i=0}^{m} \binom{k+i}{k} (1-x)^i 
	  + (1-x)^{m+1} \sum_{i=0}^{k} \binom{m+i}{m} x^i
\end{equation} 
for integers  $ k, m \ge 0 $.  Many different proofs of  (\ref{eq:CB})  are known.  See 
 \cite{KS} for a detailed account.    For the case  $ m=k $,  (\ref{eq:CB})  is frequently 
called the Daubechies identity.  See \cite{Zeilberger}.  As described in  \cite{KSII}, 
the Chaundy and Bullard inequality has roots going back three centuries.
In what follows we present some  ramifications of  (\ref{eq:CB}). 
In the first part we discuss extensions to several variables and relations with other identities.
In the second part we obtain additional identities with more parameters.

\section{The homogeneous form of (\ref{eq:CB})  }
\setcounter{section}{1}	  

 \ The homogeneous identity 
\begin{equation}							\label{eq:HOM}
x^{m+1} y^{k+1} = \sum_{i=0}^{m} \binom{k+i}{k} x^{m-i+1} \left( \frac{xy}{x+y} \right)^{k+i+1}  
	 + \  \sum_{i=0}^{k} \binom{m+i}{m} y^{k-i+1} \left( \frac{xy}{x+y} \right)^{m+i+1}  
\end{equation}
deserves attention for its own sake and has several interesting conclusions:

\noindent
(a) \  If we divide  (\ref{eq:HOM})  by  $ x^{m+1} y^{k+1} $ and choose  $ x+y = 1 $, 
we get  identity  (\ref{eq:CB}). 

\noindent
(b) \  Another conclusion of (\ref{eq:HOM}) is an identity given by  Graham, Knuth and 
Patashnik   \cite[p. 246]{Graham-Knuth-Patashnik}:
\  {\it  If $ xy=x+y  $  then }
\begin{equation}   							\label{eq:Knuth} 
	x^{m+1} y^{k+1} = \sum_{i=0}^{m} \binom{k+i}{k} x^{m-i+1} 
			+ \sum_{i=0}^{k} \binom{m+i}{m} y^{k-i+1} .
\end{equation} 
It is not difficult to see that  (\ref{eq:CB})--(\ref{eq:Knuth})  are all equivalent.

We suggest a proof of  (\ref{eq:HOM})  which conveniently generalizes to more than two variables.
Let us apply 
$ 
  \left( \!  -\dfrac{\partial} {\partial x} \right)^m  
		\mkern-9 mu 
  \left( \!  -\dfrac{\partial} {\partial y} \right)^k \
$ 
to the identity 
\begin{equation}							\label{eq:Id-1}
  		\frac{1}{x y} = \frac{1}{ x(x+y) } + \frac{1}{ y(x+y) }  . 
\end{equation}
First we apply   $ \left( -{\partial} / {\partial x} \right)^m  $.  To the term  
$  \  \dfrac{1}{ x(x+y) }  \ $  we use the Leibnitz formula 
	$  \  (fg)^{ (m) } = \sum_{i=0}^m \binom{m}{i} f^{ (m-i) } g^{ (i) }  $ :
$$
\frac{ m! } { x^{m+1} y } 
= \sum_{ i=0 }^m \binom{m}{i} \frac{ (m-i)! }{ x^{m-i+1} } \frac{ i! }{ (x+y)^{i+1} }
	+  \frac{ m! }{ y(x+y)^{m+1} } 
= m !\sum_{ i=0 }^m  \frac{ 1 }{ i! x^{m-i+1} } \frac{ i! }{ (x+y)^{i+1} }
	+  \frac{ m! }{ y(x+y)^{m+1} }  .
$$
Next, applying  $ \left(-{\partial} / {\partial y} \right)^k  $; 
\begin{equation}						\label{eq:INV} 
\begin{aligned} 
\frac{ m! \, k! } { x^{m+1} y^{k+1} } 
& = m! \sum_{ i=0 }^m  \frac{ 1 }{ i! x^{m-i+1} } \frac{ (i+k)! }{ (x+y)^{i+k+1} }
     + \sum_{ j=0 }^k \binom{k}{j} \frac{ (k-j)! }{ y^{k-j+1} }   \frac{ (m+j)! }{ (x+y)^{m+j+1} } 
											\\
& = m! \, k! \sum_{ i=0 }^m \binom{i+k}{k} \frac{ 1 }{ x^{m-i+1} (x+y)^{i+k+1} } 
    + k! \, m! \sum_{ j=0 }^k \binom{m+j}{m} \frac{ 1 }{ y^{k-j+1} (x+y)^{m+j+1} }  .
\end{aligned}
\end{equation}
Replacing  $ x, y $  by  $ x^{-1}, y^{-1} $, respectively, one gets  (\ref{eq:HOM}). 
\qed

Note that the  case  $ m=k=1 $  of  (\ref{eq:INV}),  namely 
$ \displaystyle
\  \frac{1} { x^{2} y^{2} } 
 	= \left( \frac{1}{x^2} + \frac{1}{y^2} \right) \frac{1}{ (x+y)^2 } 
 	+ \left( \frac{1}{x} + \frac{1}{y} \right) \frac{2}{ (x+y)^3 }   \ 
		\rule[-6 mm]{0 mm}{11 mm}
$, 
plays a central role in the develpment of the theory of G. Eisenstein about periodic 
functions. See  \cite[p. 252]{Ei}.


\section{A generalization for  $n$  variables }


We propose a homogeneous  identity with  $n$  variables which generalizes both 
(\ref{eq:HOM})  and  (\ref{eq:Knuth}).   Our method uses only very elementary tools of 
analysis.  We apply the differential operator  
$ \displaystyle 
  \left(-\frac{\partial}{\partial x_1} \right)^{m_1}   \cdots 
  \left(-\frac{\partial}{\partial x_n} \right)^{m_n} 
$
to the elementary identity 
\begin{equation}					\label{eq:basic-Id}
\frac{1}{ x_1 x_2 \cdots x_n } 
= \sum_{t=1}^n  \frac{1}
{ x_1 \cdots  \left\langle {{\displaystyle  x_t} \atop \hbox{\rm skipped} } \right\rangle  \cdots  x_n 
		(x_1 + \ldots + x_n) } , 
\end{equation} 
which generalizes   (\ref{eq:Id-1}).  The result of applying the operator to the left 
hand side of  (\ref{eq:basic-Id})  is  
\begin{equation}						\label{eq:left}
  \frac{ m_1 ! m_2 ! \ldots m_n ! }{ x_1^{ m_1 + 1 } x_2^{m_2 + 1} \cdots x_n^{m_n + 1} } . 
\end{equation}
On the right hand side of  (\ref{eq:basic-Id})  we differentiate each term separately, 
i.e., take a fixed  $ t $,  $  1 \le t \le n $,  and calculate 
\begin{equation}						\label{eq:right-t} 
  \left(-\frac{\partial}{\partial x_1} \right)^{m_1}   \cdots 
  \left(-\frac{\partial}{\partial x_n} \right)^{m_n} 
  \frac{1}{ x_1 \cdots  \left\langle {{\displaystyle  x_t} \atop \hbox{\rm skipped} } \right\rangle  \cdots x_n 
		(x_1 + \ldots + x_n) } . 
\end{equation}
Since the variable  $ x_t $  appears only in one factor of the denominator of  (\ref{eq:right-t}),  
while each other   $ x_j $,  $ j \neq t $  appears in two factors, we apply  
	$ \left(-{\partial} / {\partial x_t} \right)^{m_t} $ 
first and get 
\begin{equation}						\label{eq:right-diff-t} 
\frac{m_t!}{ x_1 \cdots  \left\langle {{\displaystyle  x_t} \atop \hbox{\rm skipped} } \right\rangle  \cdots x_n 
		(x_1 + \ldots + x_n)^{m_t + 1 } } . 
\end{equation}
Next we apply  $ \prod_{ j \neq t } \left(-{\partial} / {\partial x_j} \right)^{m_j} $ 
to  (\ref{eq:right-diff-t}).  By the Leibnitz formula 
we get that  (\ref{eq:right-t}) equals 
\begin{equation}						\label{eq:right-side}
\begin{aligned}
&
\prod_{ j \neq t } \left( -\frac{\partial}{\partial x_j} \right)^{m_j}
\frac{m_t!}{ x_1 \cdots  \left\langle {{\displaystyle  x_t} \atop \hbox{\rm skipped} } \right\rangle  \cdots x_n 
		(x_1 + \ldots + x_n)^{m_t + 1 } } 		\\
&  
\mkern-19 mu 
= 
\ssum_{ 0 \le i_j \le m_j \atop j \neq t }
\binom{m_1}{i_1} \frac{ (m_1 - i_1)! }{ x_1^{m_1 - i_1 + 1 } }
\cdots   \Big\langle { {\displaystyle  i_t, x_t } \above 0pt {\rm skipped} }  \Big\rangle   \cdots
\binom{m_n}{i_n} \frac{ (m_n - i_n)! }{ x_n^{m_n - i_n + 1 } }
\frac{ (m_t + i_1 + \ldots + i_{t-1} + i_{t+1} + \ldots + i_n)! }
{ (x_1 + \ldots + x_n)^{ m_t + i_1 + \ldots + i_{t-1} + i_{t+1} + \ldots + i_n + 1 } }  \ . 
\end{aligned}
\end{equation}
Summing  (\ref{eq:right-side})  for  $ t=1, \ldots, n $ and comparing with  (\ref{eq:left}),  
results  
\begin{equation}
\begin{aligned}							\label{eq:inv-id}
\frac{ 1 }{ x_1^{ m_1 + 1 } x_2^{m_2 + 1} \cdots x_n^{m_n + 1} } 
& =
\sum_{t=1}^n  
\Bigg[
\ssum_{ 0 \le i_j \le m_j \atop j \neq t }
\frac{ ( i_1 + \ldots + i_{t-1} + m_t + i_{t+1} + \ldots + i_n)! }
     { i_1!  \ldots  i_{t-1}! m_t! i_{t+1}! \ldots  i_n! }		\\
&
\times  \frac{ 1 }{ x_1^{m_1 - i_1 + 1 } 
	\ldots   \left\langle {{\displaystyle  x_t} \atop \hbox{\rm skipped} } \right\rangle   \ldots
	x_n^{m_n - i_n + 1 }  (x_1 + \ldots + x_n)^{ i_1 + \ldots + i_{t-1} + m_t + i_{t+1} + \ldots + i_n + 1 } } 
 \Bigg] . 
\end{aligned}
\end{equation}
Finally we replace   $ x_i $  by  $ x_i^{-1} $  and use for two of the basic symmetric polynomials 
in  $n$  variables the notation 
$$	 S_{n,n}(x_1, \ldots, x_n) = x_1 \cdots x_n , \quad 
	 S_{n-1,n}(x_1, \ldots, x_n) 
	   = \sum_{t=1}^n x_1  \cdots \left\langle { {\displaystyle  x_t} \atop \hbox{\rm skipped} } \right\rangle \cdots x_n .
$$
Then  (\ref{eq:inv-id})  becomes our main homogeneous identiy 
\begin{equation}						\label{eq:n-powers} 
\begin{aligned}
 x_1^{m_1 + 1} x_2^{m_2 + 1} \cdots x_n^{m_n + 1} 
 =  & 
\sum_{t=1}^n   \Bigg[ 
\ssum_{  \begin{array}{c}  0 \le i_j \le m_j 	\\   j \neq t    \end{array}  } 
\frac{ ( i_1 + \cdots + i_{t-1} + m_t +i_{t+1} + \cdots i_n )! }
     { i_1 ! \cdots  i_{t-1}! \, m_t ! \, i_{t+1}!   \cdots i_n ! } 
									\\
&  
\times \,
	x_1 ^ {m_1 - i_1 + 1}  
	\cdots \left\langle {{\displaystyle  x_t} \atop \hbox{\rm skipped} } \right\rangle  \cdots 
        x_n ^ {m_n - i_n + 1} 
\left( \frac{ S_{n,n} }{ S_{n-1,n} } \right)^{ i_1 + \cdots i_{t-1} + m_t + i_{t+1} + \cdots i_n + 1 } 
\Bigg] . 
\end{aligned}
\end{equation}

\noindent
{\bf  Examples.} \  For $ n=2 $,  (\ref{eq:n-powers})  reduces to  
\begin{equation}						\label{eq:2-powers} 
x_1^{m_1 + 1} x_2^{m_2 + 1} 
= 
\sum_{ i_2 = 0 }^{ m_2 } \binom{m_1 + i_2}{m_1} x_2^{ m_2 - i_2 + 1 } 
\left( \frac{ x_1 x_2 }{ x_1 + x_2 } \right)^{ m_1 + i_2 + 1 }
+
\sum_{ i_1 = 0 }^{ m_1 } \binom{m_2 + i_1}{m_2} x_1^{ m_1 - i_1 + 1 } 
\left( \frac{ x_1 x_2 }{ x_1 + x_2 } \right)^{ i_1 + m_2 + 1 } , 
\end{equation}
i.e.,  (\ref{eq:HOM}).

Assuming the equality  $ S_{n,n} = S_{n-1,n} $,  identity  (\ref{eq:n-powers})  
implies a  $n$-variable analogue to  (\ref{eq:Knuth}). For $ n=3 $  it is:     \\
%
%
{\it  If  \  $ xyz=xy+yz+zx  $,  \  then }
\begin{equation} 						\label{eq:Knuth3}
\begin{aligned}
 x^{m_1 + 1}  &   y^{m_2 + 1} z^{m_3 + 1} 
   = 
\sum_{ j \le m_2 , \,  k \le m_3 } 
\frac{ ( m_1 + j + k )! } { m_1! \, j! \, k! } 
	y^{ m_2 - j + 1 } z^ {m_3 - k + 1}   
								\\
& + 
\sum_{ k \le m_3 , \,  i \le m_1 } 
\frac{ ( i+ m_2 + k )! } { i! \, m_2! \, k! } 
	z^ {m_3 - k + 1} x^{ m_1 - i + 1 }  			
+ \sum_{ i \le m_1 , \,  j \le m_2 }  
\frac{ ( i +j + m_3 )! } { i ! \, j! \, m_3 ! } 
	x^ {m_1 - i + 1}   y^{ m_2 - j + 1 } 	 .
\end{aligned}
\end{equation}

If we divide  (\ref{eq:n-powers})  by  $ x_1^{m_1 + 1}  \cdots x_n^{m_n + 1} $   and 
take  $ S_{n-1,n} = 1 $,  we get a  $n$-variable analogue of the identity of 
Chaundy and Bullard.  For  $ n=3 $  it is:	\ 
%
%
{\it If  \  $ xy+yz+zx = 1 $,  \  then }
\begin{equation} 						\label{eq:S2=1}
\begin{aligned} 
(yz)^{ m_1 + 1 } 
\sum_{ j \le m_2 , \,  k \le m_3 } 
\frac{ ( m_1 + j + k )! } { m_1! \, j! \, k! } 
	y^k  z^j  x^{ j + k }	
+ 
(zx)^{ m_2 + 1 } 
\sum_{ k \le m_3 , \,  i \le m_1 } 
\frac{ ( i+ m_2 + k )! } { i! \, m_2! \, k! } 
	z^i x^k y^{i+k}		&			\\			
+ (xy)^{ m_3 + 1 } 
\sum_{ i \le m_1 , \,  j \le m_2 } 
\frac{ ( i +j + m_3 )! } { i ! \, j! \, m_3 ! } 
	x^j  y^i  z^{ i + j } 
&  = 1	.
\end{aligned}
\end{equation}
\qed

\bigskip

The change of variables 
$	u_t = x_1 \cdots  x_{t-1}  x_{t+1}  \cdots x_n $,  \ $  t=1,\ldots, n $, 
and the inverse transformation 
$$  
	x_t = \frac{ (u_1 \cdots u_n)^{  {1} / {(n-1)} } }{ u_t } ,  
			\qquad  	  t=1,\ldots, n ,		
$$ 
yield   $  \  S_{n-1,n}(x_1, \ldots, x_n) = u_1 + \ldots + u_n $,
	$  \  S_{n,n}(x_1, \ldots, x_n) = ( u_1 \ldots u_n )^{ 1/(n-1) } $. 
After some elementary calculations this transforms identity  (\ref{eq:n-powers})  into 
\begin{equation} 							\label{eq:transformed}
\begin{aligned}
 {  (u_1 + \cdots + u_n) }^{ m_1 + \cdots + m_n + 1 }  
  =  &
\sum_{t=1}^n   \Biggl[  u_t^{ m_t + 1 } \!\!\!
\ssum_{  \begin{array}{c}   0 \le i_j \le m_j 	\\   j \neq t    \end{array}  } 
\frac{ ( i_1 + \cdots + i_{t-1} + m_t +i_{t+1} + \cdots i_n )! }
     { i_1 ! \cdots  i_{t-1}! \, m_t ! \, i_{t+1}!   \cdots i_n ! } 
									\\
& 
\times	  
{u_1}^{i_1} \cdots \left\langle {{\displaystyle  u_t} \atop \hbox{\rm skipped} } \right\rangle  \cdots  {u_n}^{i_n} \, 
{  (u_1 + \cdots + u_n) }^{ \sum_{j \neq t } (m_j - i_j) }
\Biggr] .
\end{aligned}
\end{equation}

(\ref{eq:transformed})  is precisely equation  (10.2)  of  \cite{KS}.  Two proofs 
of this result are given in  \cite{SIAM},  one by a probabilistic argument and the 
other by using generating functions.


\section{Another generalization of CB}

The next identity is another generalization of  (\ref{eq:CB})  which depends on three 
independent parameters:


\noindent
{\it  Let $ m - r + k - \ell = 0 $,   $ m, r, k, \ell $ positive integers. Then }
\begin{equation}						\label{eq:3-parms}
(1-x)^{r+1} \sum_{i=0}^{k} \binom{m+i}{r} x^{i+m-r} 
\ + \ x^{\ell+1} \sum_{i=0}^{m} \binom{k+i}{ \ell } (1-x)^{i+k-\ell} 
= \left\{  \begin{array}{ll} 
	\displaystyle
	 1   -  \!\!\! \sum_{i=0}^{m-r-1} \! \binom{m}{i} x^i (1-x)^{m-i}  
						&  \!\!   \mbox{if \, $ m - r > 0 $, }  \\
	\displaystyle
	 1  -  \!\!  \sum_{i=0}^{ k - \ell - 1 }  \binom{k}{i} (1-x)^{i} x^{k-i} 	
						&  \!\!   \mbox{if \ $  k - \ell > 0 $, }
	 \rule{0 mm}{8 mm}  							\\
	1					&   \!\! \mbox{if \, $ m = r , k=l $ .} 
	\rule{0 mm}{6 mm} 
\end{array}   \right .  
\end{equation}

In the previous version of this manuscript  (\ref{eq:3-parms})  was proved by elementary 
methods using ideas presented in  \cite{AE}.  Professor T. Koornwinder kindly
brought to our attention a shorter proof of  (\ref{eq:3-parms}), which follows hereby:

Let us verify the case  $ m-r = \ell - k > 0 $.  In the second sum on the left hand side 
the terms are nonzero only when  $ k+i \ge \ell $,  hence it is sufficient to sum only for 
$ i \ge  \ell - k = m-r $.  We change the summation index in the first sum on the left 
to  $ j = i + (m - r) $  and in the second sum to   $ j = i - (m - r) $.  By repeated use
of  $ m-r = \ell - k $,  the left side becomes 
$$
(1-x)^{r+1} \sum_{ j = m - r }^{ k + m - r } \binom{ j + r }{r} x^{j} 
\ + \ 
 x^{\ell+1} \sum_{ j = 0 }^{ r } \binom{ j + \ell }{ \ell } (1-x)^{j} . 
$$
Let us rewrite this as 
$$
(1-x)^{r+1} 
\left[ 
	\sum_{ j = 0 }^{ k + m - r } \binom{ j + r }{r} x^{j} 
      -
	\sum_{ j = 0 }^{ m - r - 1 } \binom{ j + r }{r} x^{j}
\right]
\ + \ 
 x^{\ell+1} \sum_{ j = 0 }^{ r } \binom{ j + \ell }{ \ell } (1-x)^{j} , 
$$
and rearrange it to 
$$
\left[
(1-x)^{r+1}  \sum_{ j = 0 }^{ \ell } \binom{ j + r }{r} x^{j} 
\ + \ 
 x^{\ell+1} \sum_{ j = 0 }^{ r } \binom{ j + \ell }{ \ell } (1-x)^{j} 
\right]
-
(1-x)^{r+1}  \sum_{ j = 0 }^{ m - r - 1 } \binom{ j + r }{r} x^{j} . 
$$
The first two sums total to 1 by the original Chaundy-Bullard identity, so  (\ref{eq:3-parms}) 
will follow if one shows that
$$
1 - (1-x)^{r+1}  \sum_{ j = 0 }^{ m - r - 1 } \binom{ j + r }{r} x^{j} 
=
1 -  \sum_{i=0}^{m-r-1} \! \binom{m}{i} x^i (1-x)^{m-i} ,
$$
i.e.,
$$
1 - (1-x)^{r+1}  \sum_{ j = 0 }^{ m - r - 1 } \binom{ j + r }{r} x^{j} 
=
1 - (1-x)^{r+1}  \sum_{i=0}^{m-r-1} \! \binom{m}{i} x^i (1-x)^{m-r-i-1} ,
$$
But the remaining 
$$
  \sum_{ j = 0 }^{ m - r - 1 } \binom{ j + r }{r} x^{j} 
=
  \sum_{i=0}^{m-r-1} \! \binom{m}{i} x^i (1-x)^{m-r-i-1} 
$$
is precisely equation (2.7) of  \cite{KS},  hence  (\ref{eq:3-parms})  follows.
\qed


\bigskip
\noindent
{\bf Acknowkedgment.}  The authors would like to thank Professor T. Koornwinder 
for the stimulating correspondence.



\bigskip

\noindent\textit
{Department of Mathematics, Technion --- I.I.T., Haifa 32000, Israel	\\ 
 dova@tx.technion.ac.il		\\	 elias@tx.technion.ac.il}

\end{document}